\documentclass[10pt]{article}
\usepackage{amssymb}
\usepackage{fullpage}
\usepackage{amsthm}
\usepackage{amsmath}
\usepackage{mathrsfs}
\usepackage{natbib}

\begin{document}
\newcommand{\beq}{\begin{eqnarray}}
\newcommand{\eeq}{\end{eqnarray}}
\newcommand{\beas}{\begin{eqnarray*}}
\newcommand{\enas}{\end{eqnarray*}}
\newcommand{\bea}{\begin{eqnarray}}
\newcommand{\ena}{\end{eqnarray}}
\newcommand{\nn}{\nonumber}
\newcommand{\ignore}[1]{}
\newcommand{\qmq}[1]{\quad \mbox{#1} \quad}
\newtheorem{theorem}{Theorem}[section]
\newtheorem{corollary}{Corollary}[section]
\newtheorem{conjecture}{Conjecture}[section]
\newtheorem{proposition}{Proposition}[section]
\newtheorem{remark}{Remark}[section]
\newtheorem{lemma}{Lemma}[section]
\newtheorem{definition}{Definition}[section]
\newtheorem{condition}{Condition}[section]
\newcommand{\pf}{\noindent {\bf Proof:} }
\def\blfootnote{\xdef\@thefnmark{}\@footnotetext}
\title{Concentration of measure for the number of isolated vertices in the Erd\H{o}s-R\'{e}nyi random graph by size bias couplings}
\author{Subhankar Ghosh\footnote{Department of Mathematics, University of Southern California, Los Angeles, CA 90089, USA \texttt{subhankg@usc.edu}}, Larry Goldstein\footnote{Department of Mathematics, University of Southern California, Los Angeles, CA 90089, USA \texttt{larry@math.usc.edu}}\,\,\thanks{research partially supported by NSA grant H98230-11-1-0162}
 and Martin Rai\v{c}\footnote{University of Ljubljana and University of Primorska, Slovenia \texttt{martin.raic@fmf.uni-lj.si}}}

\date{}
\maketitle
\begin{abstract}
A concentration of measure result is proved for the number of isolated vertices $Y$ in the Erd\H{o}s-R\'{e}nyi random graph model on $n$ edges with edge probability $p$. When $\mu$ and $\sigma^2$ denote the mean and variance of $Y$ respectively, $P((Y-\mu)/\sigma\ge t)$ admits a bound of the form $e^{-kt^2}$ for some constant positive $k$ under the assumption $p \in (0,1)$ and $np\rightarrow c \in (0,\infty)$ as $n \rightarrow \infty$.
The left tail inequality
\beas
P\left(\frac{Y-\mu}{\sigma}\le -t\right)&\le& \exp\left(-\frac{t^2\sigma^2}{4\mu}\right)
\enas
holds for all $n \in \{2,3,\ldots\},p \in (0,1)$ and $t \ge 0$. The results are shown by coupling $Y$ to a random variable $Y^s$ having the $Y$-size biased distribution, that is, the distribution characterized by $E[Yf(Y)]=\mu E[f(Y^s)] $ for all functions $f$ for which these expectations exist.
\end{abstract}
\long\def\symbolfootnote[#1]#2{\begingroup%
\def\thefootnote{\fnsymbol{footnote}}\footnote[#1]{#2}\endgroup}
\symbolfootnote[0]{2000 {\em Mathematics Subject Classification}: Primary
60E15\ignore{Inequalities; stochastic orderings};
Secondary 60C05\ignore{Combinatorial probability},60D05\ignore{Geometric probability and stochastic geometry}.}
\symbolfootnote[0]{{\em Keywords}: Large deviations; graph degree; size biased couplings}

\section{Introduction and main result}

For some $n \in \{1,2,\ldots,\}$ and $p \in (0,1)$ let $K$ be the Erd\H{o}s-R\'{e}nyi random graph on the vertices ${\cal V}=\{1,2,\ldots,n\}$ and edge success probability $p$, that is, with edge indicators $X_{uv}=1(u,v\in\mathcal{V}:\mbox{$\overline{uv}$ is an edge in $K$})$ independent random variables with the Bernoulli($p$) distribution for all $u\neq v$. We set $X_{vv}=0$
for all $v \in {\cal V}$. Recall that the degree of a vertex $v \in {\cal V}$, denoted by $d(v)$, is the number of edges incident on $v$. Hence,
\beas 
d(v)=\sum_{w \in {\cal V}} X_{vw}.
\enas
Many authors have studied the distribution of
\bea
Y=\sum_{v \in {\cal V}} 1(d(v)=d)\label{Ysum}
\ena
counting the number of vertices $v$ of $K$ having some fixed degree $d$. We derive upper bounds, for fixed $n$ and $p$, on the tail probabilities of the number of isolated vertices of $K$, that is, for $Y$ in (\ref{Ysum}) for the case $d=0$, which counts the number of vertices having no incident edges.

For $d$ in general, and $p$ depending on $n$, previously in \cite{Karonski87},
the asymptotic normality of $Y$ was shown when $n^{(d+1)/d} p\rightarrow \infty$ and $n p\rightarrow 0$, or $n p \rightarrow \infty$ and $n p-\log n
-d \log \log n \rightarrow - \infty$; see also \cite{Palka84} and \cite{Bollobas85}. For the case $d=0$ of isolated vertices, \cite{barbour82} and \cite{barbour} show that $Y$ is asymptotic normal if and only if $n^2p \rightarrow \infty$ and $np - \log n \rightarrow -\infty$.

Here we study the distribution of $Y$ using a size bias coupling that was used in \cite{golrinott} to study the rate of convergence to the multivariate normal distribution for a vector whose components count the number of vertices of some fixed degrees. In \cite{kordeki}, the mean $\mu$ and variance $\sigma^2$ of $Y$ for the particular case $d=0$ are computed as
\bea \label{graph:meanvar}
\mu=n(1-p)^{n-1},\quad\mbox{and}\quad \sigma^2=n(1-p)^{n-1}(1+np(1-p)^{n-2}-(1-p)^{n-2}) \qmq{for $n \ge 2$.}
\ena
In the same paper, Kolmogorov distance bounds to the normal of order $O(\mbox{Var}(Y)^{-1/2})$
were obtained.

\cite{oc} showed that an asymptotic large deviation principle holds for $Y$. \cite{raic} obtains nonuniform large deviation bounds for mean zero, variance one random variables in some generality, and applies his results to the case of counting the number of isolated vertices with $W=(Y-\mu)/\sigma$, yielding the bound
\bea \label{raic3}
\frac{P(W \ge t)}{1-\Phi(t)}\le e^{t^3\beta(t)/6}(1+Q(t)\beta(t)) \quad \mbox{for all $t \ge 0$,}
\ena
where $\Phi(t)$ denotes the distribution function of a standard normal variate,
\beas
Q(t)=\frac{12}{\sqrt{2\pi}}+\frac{23}{2}t+\frac{11\sqrt{2\pi}}{2}t^2,
\enas
and
\beas
\beta(t)=\frac{n}{6\sigma^3}(13+43np+27(np)^2)
\exp\left\{\frac{(8+4np)t}{\sigma}+2np(e^{t/\sigma}-1)\right\}.
\enas
Still from \cite{raic}, when $np \rightarrow c$ as $n \rightarrow \infty$, (\ref{raic3}) holds for all $n$ sufficiently large with
\bea \label{raic.beta}
\beta(t) = \frac{C_1}{\sqrt{n}}\exp\left( \frac{C_2 t}{\sqrt{n}} + C_3 (e^{C_4t/\sqrt{n}}-1)\right)
\ena
for some unspecified constants $C_1,C_2,C_3$ and $C_4$ depending only on $c$. For $t$ of order $\sqrt{n}$, for instance, the function $\beta(t)$ will be of order $1/\sqrt{n}$ as $n \rightarrow \infty$, allowing an asymptotic approximation of the deviation probability $P(W \ge t)$ by the normal, to within some factors.

In Theorem \ref{thm-graph-iso} we supply a bound that likewise holds also for all $n$, and that also gives somewhat more explicit information on the rate of tail decay. In particular, we see from (\ref{left-graph}) that the standardized variable $W$ has a left tail that is bounded above by $\exp\{-t^2 \sigma^2/(4\mu)\}$. Moreover, the right tail also exhibits similar bounds over some parameter regions, with a worst case bound there of order $\exp\{-\rho t\}$ by (\ref{right-graph-1}) for some $\rho>0$, but see also Corollary \ref{cor:ksquared} for a further improvement in the regime where $np$ converges to a nonzero constant.

\begin{theorem}\label{thm-graph-iso}
For $n \in \{2,3,\ldots\}$ and $p \in (0,1)$ let $K$ denote the random graph on $n$ vertices where each edge is present with probability $p$, independently of all other edges, and let $Y$ denote the number of isolated vertices in $K$, having mean $\mu$ and variance $\sigma^2$, as given in (\ref{graph:meanvar}). Let $M(\theta)=E\exp(\theta(Y-\mu)/\sigma)$ be the moment generating function of the standardized $Y$ variable. Then, letting
\beas 
\gamma_s= e^s(pe^s+1-p)^{n-2}(npe^s+1-p)+(n-1)p+1 \qmq{and} H(\theta) = \frac{\mu}{2\sigma^2} \int_0^\theta s \gamma_{s/\sigma} ds
\enas
we have $M(\theta) \le \exp H(\theta)$ for all $\theta \ge 0$, and for all $t>0$,
\bea
\label{right-graph}
P\left(\frac{Y-\mu}{\sigma}\ge t\right) \le \inf_{\theta \ge 0} \exp(-\theta t+H(\theta)).
\ena

For all $\theta \le 0$ we have $M(\theta) \le \exp(\mu \theta^2/\sigma^2)$, and for all $t>0$,
\bea \label{left-graph}
P\left(\frac{Y-\mu}{\sigma}\le -t\right)&\le& \exp\left(-\frac{t^2\sigma^2}{4\mu}\right).
\ena
\end{theorem}

Though integration shows that we may explicitly write
\beas
H(\theta)= \frac{\mu}{2\sigma^2}\left(\frac{n p  ((n-1) p+1)\theta ^2 + 2 \sigma ^2  -2 \sigma  \left(p e^{\theta /\sigma } (\sigma -n \theta )
+\sigma(1-p)
\right) \left(p \left(e^{\theta /\sigma }-1\right)+1\right)^{n-1}}{2np}\right),
\enas
the integral formula for $H(\theta)$ in the theorem appears simpler to handle.


Useful bounds for the minimization in (\ref{right-graph})  may be obtained by restricting
to $\theta \in [0,\theta_0]$ for some $\theta_0$. In this case, as $\gamma_{s/\sigma}$ is an increasing function of $s$,
we have
\beas
H(\theta) \le  \frac{\mu}{4\sigma^2} \gamma_{\theta_0/\sigma}\theta^2  \quad \mbox{for $\theta \in [0,\theta_0]$.}
\enas
The quadratic $-\theta t + \mu \gamma_{\theta_0/\sigma} \theta^2/(4\sigma^2)$ in $\theta$ is minimized at $\theta=2t\sigma^2/(\mu \gamma_{\theta_0/\sigma})$. When this value falls in $[0,\theta_0]$ we obtain the first bound in (\ref{right-graph-1}), and setting $\theta=\theta_0$ yields the second, thus,
\bea \label{right-graph-1}
P\left(\frac{Y-\mu}{\sigma}\ge t\right) \le \left\{
\begin{array}{ll}
\exp(-\frac{t^2 \sigma^2}{\mu \gamma_{\theta_0/\sigma} }) & \mbox{for $t \in [0,\theta_0 \mu \gamma_{\theta_0/\sigma} /(2 \sigma^2)]$}\\
\exp(-\theta_0 t+\frac{\mu \gamma_{\theta_0/\sigma} \theta_0^2}{4 \sigma^2}) & \mbox{for $t \in (\theta_0 \mu \gamma_{\theta_0/\sigma} /(2 \sigma^2),\infty)$.}
\end{array}
\right.
\ena

Inequality (\ref{right-graph-1}) and the boundedness of $Y$ yields the following useful corollary.
\begin{corollary}
\label{cor:ksquared}
For all $c \in (0,\infty)$ there exists a positive constant $k$ depending only on $c$ such that when $p \in (0,1)$ and $np \rightarrow c \in (0,\infty)$ as $n \rightarrow \infty$,
\beas
P((Y-\mu)/\sigma \ge t) \le \exp(-kt^2)
\enas
for all $t \ge 0$ and all $n \ge 2$.
\end{corollary}

\proof Since $Y$ can be no more than $n$, and $\sigma^2$ increases at rate $n$ when $np \rightarrow c$, there exists a positive constant $a_0$ such that
\beas
\frac{Y-\mu}{\sigma} \le \frac{n}{\sigma} \le a_0 \sqrt{n}.
\enas
Hence $P((Y-\mu)/\sigma) \ge t)=0$ for all $t > a_0 \sqrt{n}$.

For any given $n$ let $\theta_n=a_0\sqrt{n}\sigma^2/\mu$. Then, as $\gamma_s \ge 2$ for all $s \ge 0$, we have $(\theta_n \mu \gamma_{\theta_0/\sigma})/(2 \sigma^2) \ge a_0 \sqrt{n}$, so the first bound in (\ref{right-graph-1}) applies for all $t \le a_0 \sqrt{n}$. Note that $\theta_n/\sigma = a_0\sqrt{n} \sigma/\mu$ converges to a positive constant, implying the convergence of $\gamma_{\theta_n/\sigma}$, and hence that of $\sigma^2/(\mu \gamma_{\theta_n/\sigma})$, also to a positive constant. Since $\sigma^2/(\mu \gamma_{\theta_n/\sigma})$ is positive for all $n \ge 2$, we see that the claim of the corollary holds for all $k$ in the nonempty interval $(0,\inf_n \sigma^2/(\mu \gamma_{\theta_n/\sigma}))$. \qed

In the asymptotic of Corollary \ref{cor:ksquared}, for, say $t=a \sqrt{n}$, the function $\beta(t)$ of (\ref{raic.beta}) behaves like $C/\sqrt{n}$, so
the bound (\ref{raic3}) also gives useful information for some range of positive values of $a$ up to some upper limit. However as $\exp(t^3\beta(t)/6)$ behaves like $\exp(Cna^3/6)$, when multiplied by $1-\Phi(t)$, of exponential order $\exp(-a^2n/2)$, the product tends to infinity for all sufficiently large $a$, so the bound in (\ref{raic3}) may explode before the right tail of $W$ vanishes.

%

The main tool used in proving Theorem \ref{thm-graph-iso} is size bias coupling, that is, the construction of $Y$ and $Y^s$ on the same space where $Y^s$ has the $Y$-size biased distribution characterized by
\bea
E[Yf(Y)]=\mu E [f(Y^s)]\label{EWfWchar}
\ena
for all $f$ for which the expectations above exist.
In \cite{cnm} and \cite{cnm.bd}, size bias couplings were used to prove concentration of measure inequalities when $|Y^s-Y|$ can be almost surely bounded by a constant independent of the problem size. Here, in contrast, we apply the coupling for the number of isolated vertices of $K$ from \cite{golrinott}, which violates the boundedness condition. Unlike the theorem used in \cite{cnm} and \cite{cnm.bd}, which can be applied to a wide variety of situations under a bounded coupling assumption, it seems that cases where the coupling is unbounded, such as the one we consider here, need
application specific treatment, and cannot be handled by one single general result.

Having its roots in the work of \cite{baldi}, a general prescription for constructing a variable with the size bias distribution of a sum of nonnegative variables is given in \cite{golrinott}. Helped by the fact that size biasing a nontrivial indicator random variable simply sets its value to one, specializing to nontrivial exchangeable indicators yields the following simplification as in Lemma 3.3 of \cite{golpen}.
\begin{proposition} \label{prop-golrinott}
Suppose $Y=\sum_{v \in {\cal V}}X_v$, a finite sum of nontrivial exchangeable
Bernoulli variables $\{X_v, v \in {\cal V}\}$, and that for $w \in {\cal V}$
the variables
$\{X_v^w, v \in {\cal V}\}$ have joint distribution
\beas 
{\cal L}(X_v^w, v \in {\cal V})={\cal L}(X_v, v \in {\cal V}|X_w=1).
\enas
Then
$$
Y^w=\sum_{v \in {\cal V}} X_v^w
$$
has the $Y$ size biased distribution $Y^s$, as does the mixture $Y^V$ when $V$ is a random index with values in ${\cal V}$, chosen independent of all other variables.
\end{proposition}

Construction of the variable $Y^s$ is not enough for our purposes; one must couple $Y^s$ to $Y$. However, Proposition \ref{prop-golrinott} suggests a natural coupling. Given the exchangeable indicators $\{X_v,v \in {\cal V}\}$ that sum to $Y$, choose a summand uniformly and independently. If the summand value is already one, set $Y^s=Y$. Otherwise, set this variable to one, and `adjust' the remaining variables to have their conditional distribution given that this variable takes on the value one. By Proposition \ref{prop-golrinott} the
sum $Y^s$ of these new variables has the $Y$-size biased distribution.

\section{Proof of Theorem 1.1}
For any graph with vertex set ${\cal V}$, for $v \in {\cal V}$ we let $N(v)$ denote the set of neighbors of $v$,
\beas
N(v)=\{w \in {\cal V}: X_{vw}=1\},
\enas
where $X_{vw}$ is the indicator that there exists and edge connecting vertices $v$ and $w$. We now present the proof of Theorem \ref{thm-graph-iso}.

\proof Following Proposition \ref{prop-golrinott}, we first construct a coupling of $Y^s$, having the $Y$-size bias distribution, to $Y$. Let $K$ be given, and let $Y$ be the number of isolated vertices in $K$. From (\ref{Ysum}) with $d=0$ we see that $Y$ is the sum of exchangeable indicators. Let $V$ be uniformly chosen from ${\cal V}$, independent of the remaining variables. If $V$ is already isolated, do nothing and set $K^s=K$. Otherwise, let $K^s$ be the graph obtained by deleting all the edges connected to $V$ in $K$. By Proposition \ref{prop-golrinott}, the variable $Y^s$ counting the number of isolated vertices of $K^s$ has the $Y$-size biased distribution.

Since all edges incident to the chosen $V$ are removed in order to form $K^s$,
any neighbor of $V$ which had degree
one thus becomes isolated, and $V$ also becomes isolated if it was not so earlier.
 As $1(d(w)=0)$ is unchanged for all $w \not \in \{V\} \cup N(V)$, we have
\bea\label{diff-graph}
 Y^s-Y= d_1(V)+ 1(d(V) \not = 0)\quad \mbox{where for any $v \in {\cal V}$ we let} \quad d_1(v)=\sum_{w \in N(v)}1(d(w)=1).
\ena
In particular the coupling is monotone, that is, $Y^s\ge Y$. Further, since $d_1(V)\le d(V)$, by (\ref{diff-graph}) we have
\bea\label{bd-d1-d}
0\le Y^s-Y\le d(V)+1.
\ena

Now note that, for real $x\neq y$, the convexity of the exponential function implies
\beas
\frac{e^y-e^x}{y-x}=\int_0^1e^{ty+(1-t)x}dt \le \int_0^1(te^y+(1-t)e^x)dt =\frac{e^y+e^x}{2},
\enas
and therefore, for all real $x,y$,
\bea\label{ineq-main}
|e^x-e^y|\le |x-y|\frac{e^y+e^x}{2}.
\ena
Let $\theta \ge 0$. Using (\ref{ineq-main}) and (\ref{bd-d1-d}),
we have
\bea
E(e^{\theta Y^s}-e^{\theta Y}) &\le&  \frac{\theta}{2} E\left((Y^s-Y)(e^{\theta Y^s}+e^{\theta Y})\right)\nn \\
&=& \frac{\theta}{2}E\left(e^{\theta Y}(Y^s-Y)(e^{\theta(Y^s-Y)}+1)\right)
\le
\frac{\theta}{2}E\left(e^{\theta Y}(d(V)+1)(e^{\theta(d(V)+1)}+1)\right).\label{bd-cond}
\ena
Clearly the number of isolated vertices $Y$ is a nonincreasing function of the edge indicators $X_{vw}$, while $d(V)+1$  is a nondecreasing function of these same indicators. Hence $Y$ and $d(V)+1$ have negative correlations, that is, by the inequality of \cite{harris},
\bea \label{harris}
E[f(Y)g(d(V)+1)]\le E[f(Y)]E[g(d(V)+1)]
\ena
for any two nondecreasing real functions $f$ and $g$. In particular, when $f(x)=e^{\theta x}$ and $g(x)=x(e^{\theta x}+1)$ with $x \in [0,\infty)$, by (\ref{bd-cond}) and (\ref{harris}) we obtain
\bea
E(e^{\theta Y^s}-e^{\theta Y})&\le& \frac{\theta}{2}Ee^{\theta Y} E\left((d(V)+1)(e^{\theta(d(V)+1)}+1)\right)\nn \\
 &=& \frac{\theta}{2}Ee^{\theta Y} \left( e^\theta E\left(d(V)e^{\theta d(V)}+e^{\theta d(V)}\right)+E(d(V))+1 \right).\label{YsmYbd}
\ena

To handle the terms in (\ref{YsmYbd}), note that for any vertex $v$ the degree $d(v)$ has the Binomial($n-1,p$) distribution, and in particular
\beas
E(d(v))=(n-1)p\quad\mbox{and}\quad E(e^{\theta d(v)})= \alpha_\theta \qmq{where} \alpha_\theta=(pe^\theta+1-p)^{n-1}.
\enas
Hence, as $V$ is chosen uniformly over the vertices $v$ of $K$,
\bea
E(d(V))=(n-1)p \quad\mbox{and}\quad E(e^{\theta d(V)})=\alpha_\theta,\label{momemts-dv}
\ena
and now differentiation under the second expectation above, allowed since $d(V)$ is bounded, yields
\bea
E(d(V)e^{\theta d(V)})=\phi_\theta \qmq{where} \phi_\theta=(n-1)pe^\theta(pe^\theta+1-p)^{n-2}.\label{der-dv}
\ena
Substituting (\ref{momemts-dv})  and (\ref{der-dv}) into (\ref{YsmYbd}) yields, for all $\theta \ge 0$,
\bea \label{diff-ws-w}
E(e^{\theta Y^s}-e^{\theta Y})\le  \frac{\theta \gamma_\theta}{2}E(e^{\theta Y}) \quad \mbox{where} \quad \gamma_\theta =e^\theta (\phi_\theta+\alpha_\theta) + (n-1)p +1.
\ena

Letting $m(\theta)=E(e^{\theta Y})$, using that $Y$ is bounded to differentiate under the expectation, along with
(\ref{EWfWchar}) and (\ref{diff-ws-w}), we obtain
\bea\label{bd-graph-K}
m'(\theta)=E(Ye^{\theta Y})=\mu E(e^{\theta Y^s})\le \mu\left(1+\frac{\theta \gamma_\theta}{2}\right)m(\theta).
\ena
Standardizing $Y$, we set
\bea \label{def.M}
M(\theta)=E(\exp(\theta(Y-\mu)/\sigma))=e^{-\theta\mu/\sigma}m(\theta/\sigma),
\ena
and now by differentiating and applying (\ref{bd-graph-K}), we obtain
\beas
M'(\theta)&=&\frac{1}{\sigma}e^{-\theta\mu/\sigma}m'(\theta/\sigma)-\frac{\mu}{\sigma}e^{-\theta\mu/\sigma}m(\theta/\sigma)\nonumber\\
&\le &\frac{\mu}{\sigma}e^{-\theta\mu/\sigma}\left(1+\frac{\theta \gamma_{\theta/\sigma}}{2\sigma}\right)m(\theta/\sigma)-\frac{\mu}{\sigma}e^{-\theta\mu/\sigma}m(\theta/\sigma)\nn\\
&=& e^{-\theta\mu/\sigma}\frac{\mu \theta \gamma_{\theta/\sigma}}{2\sigma^2}m(\theta/\sigma)=\frac{\mu \theta \gamma_{\theta/\sigma} }{2\sigma^2}M(\theta). 
\enas
Since $M(0)=1$, integrating $M'(s)/M(s)$ over $[0,\theta]$ yields the bound
\beas
\log(M(\theta))\le H(\theta), \quad\mbox{or that}\quad M(\theta)\le \exp(H(\theta)) \quad \mbox{where} \quad H(\theta)= \frac{\mu}{2\sigma^2} \int_0^\theta s \gamma_{s/\sigma} ds,
\enas
proving the claim on $M(\theta)$ for $\theta \ge 0$.  Moreover, for $\theta$ nonnegative,
\beas
P\left(\frac{Y-\mu}{\sigma}\ge t\right)&\le& P\left(\exp\left(\frac{\theta(Y-\mu)}{\sigma}\right)\ge e^{\theta t}\right)
\le e^{-\theta t} M(\theta)\le \exp(-\theta t+ H(\theta)).
\enas
As the inequality holds for all $\theta \ge 0$, it holds for the infimum over $\theta \ge 0$, proving (\ref{right-graph}).

To demonstrate the left tail bound let $\theta<0$. Since $Y^s\ge Y$  and $\theta<0$, using (\ref{ineq-main}), (\ref{bd-d1-d}) and that $Y$ is a function of $K$ we obtain
\bea
E(e^{\theta Y}-e^{\theta Y^s}) \le \frac{|\theta|}{2}E\left((e^{\theta Y}+e^{\theta Y^s})(Y^s-Y)\right)
\le |\theta|E(e^{\theta Y}(Y^s-Y))
=|\theta|E(e^{\theta Y}E(Y^s-Y|K)). \label{left-tail}
\ena
By (\ref{diff-graph}) have
\bea
E(Y^s-Y|K)=\frac{1}{n}\sum_{v\in {\cal V}} (d_1(v)+1(d(v) \not =0)) \le \frac{1}{n}\sum_{v\in {\cal V}} d_1(v)+1,
\label{cond-exp}
\ena
and noting that
\beas
\sum_{v \in {\cal V}} d_1(v) = \sum_{v \in {\cal V}} \sum_{w \in N(v)} 1(d(w)=1)
= \sum_{w \in {\cal V}} \sum_{v \in N(w)} 1(d(w)=1) =
\sum_{w \in {\cal V}}|N(w)| 1(d(w)=1) =
\sum_{w \in {\cal V}}1(d(w)=1),
\enas
the number of degree one vertices in $K$, by (\ref{cond-exp}) we find that $E(Y^s-Y|K) \le 2$.

Now, by (\ref{left-tail}) and (\ref{cond-exp}),
\beas
E(e^{\theta Y}-e^{\theta Y^s})&\le& 2|\theta|E(e^{\theta Y})
\enas
and therefore, justifying differentiating under the expectation as before, applying (\ref{EWfWchar}) yields
\beas
m'(\theta)=  E(Ye^{\theta Y})= \mu E(e^{\theta Y^s})\ge \mu \left(1+2\theta \right) m(\theta).
\enas
Again with $M(\theta)$ as in (\ref{def.M}),
\beas M'(\theta)&=&\frac{1}{\sigma}e^{-\theta\mu/\sigma}m'(\theta/\sigma)-\frac{\mu}{\sigma}e^{-\theta\mu/\sigma}m(\theta/\sigma)\\
&\ge& \frac{\mu}{\sigma}e^{-\theta \mu/\sigma}((1+2\theta/\sigma)m(\theta/\sigma))-\frac{\mu}{\sigma}e^{-\theta\mu/\sigma}m(\theta/\sigma)\\
&=&\frac{2\mu \theta}{\sigma^2}M(\theta).
\enas
Dividing by $M(\theta)$, integrating over $[\theta,0]$ and exponentiating yields
\bea
\label{bd-growth-neg-log}
M(\theta) \le \exp\left(\frac{\mu \theta^2}{\sigma^2}\right),
\ena
showing the claimed bound on $M(\theta)$ for $\theta<0$.
The inequality in (\ref{bd-growth-neg-log}) implies that for all $t>0$ and $\theta<0$,
\beas
P\left(\frac{Y-\mu}{\sigma}\le -t\right)\le \exp\left(\theta t + \frac{\mu   \theta^2}{\sigma^2}\right).
\enas
Taking $\theta=-t\sigma^2/(2\mu)$ we obtain (\ref{left-graph}).\qed


\begin{thebibliography}{99}

\bibitem[Baldi et~al.(1989)Baldi, Rinott, and Stein]{baldi}
P.~Baldi, Y.~Rinott, and C.~Stein.
\newblock A normal approximations for the number of local maxima of a random
  function on a graph.
\newblock In Anderson {T.W.}, Athreya {K.B.}, and Iglehart {D.L.}, editors,
  \emph{Probability, Statistics and Mathematics, Papers in Honor of Samuel
  Karlin}, 59--81. 1989.

\bibitem[Barbour(1982)]{barbour82}
A.D.~Barbour.
\newblock Poisson convergence and random graphs.
\newblock \emph{Math. Proc. Cambridge Philos. Soc.}, 92:\penalty0 349--359, 1982.


\bibitem[Barbour et~al.(1989)Barbour, Karo\'{n}ski, and Ruci\'{n}ski]{barbour}
A.D. Barbour, M.~Karo\'{n}ski, and A.~Ruci\'{n}ski.
\newblock A central limit theorem for decomposable random variables with
  applications to random graphs.
\newblock \emph{J. Combinatorial Theory B}, 47:\penalty0 125--145, 1989.

\bibitem[Bollob{\'a}s(1985)]{Bollobas85}
B.~Bollob{\'a}s.
\newblock \emph{Random graphs}.
\newblock Academic Press Inc., London, 1985.


\bibitem[Ghosh and Goldstein(2011{\natexlab{a}})]{cnm}
S.~Ghosh and L.~Goldstein.
\newblock Concentration of measures via size biased couplings.
\newblock \emph{Probab. Th. Rel. Fields}, 149:\penalty0 271--278,
  2011{\natexlab{a}}.

\bibitem[Ghosh and Goldstein(2011{\natexlab{b}})]{cnm.bd}
S.~Ghosh and L.~Goldstein.
\newblock Applications of size biased couplings for concentration of measures.
\newblock \emph{Electronic Communications in Probability}, 16:\penalty0 70--83,
  2011{\natexlab{b}}.

\bibitem[Goldstein and Penrose(2010)]{golpen}
L.~Goldstein and M.~Penrose
\newblock Normal approximation for coverage models over binomial point processes.
\newblock\emph{Ann. Appl. Probab.}, 20:\penalty0 696-721.
\bibitem[Goldstein and Rinott(1996)]{golrinott}
L.~Goldstein and Y.~Rinott.
\newblock Multivariate normal approximations by Stein's method and size bias
  couplings.
\newblock \emph{Journal of Applied Probability}, 33:\penalty0 1--17, 1996.

\bibitem[Harris(1960)]{harris}
T.~Harris.
\newblock A lower bound for the critical probability in a certain
             percolation process.
\newblock \emph{Proceedings of the Cambridge Philosophical Society} 56:\penalty0 13�-20, 1960.

\bibitem[Karo{\'n}ski and Ruci{\'n}ski(1987)]{Karonski87}
M.~Karo{\'n}ski and A.~Ruci{\'n}ski.
\newblock Poisson convergence and semi-induced properties of random graphs.
\newblock \emph{Math. Proc. Cambridge Philos. Soc.}, 101:\penalty0 291--300,
  1987.

\bibitem[Kordecki(1990)]{kordeki}
W.~Kordecki.
\newblock Normal approximation and isolated vertices in random graphs.
\newblock In Karo\'{n}ski {M.}, Jaworski {J.}, and Ruci\'{n}ski {A.}, editors,
  \emph{Random Graphs '87}, pages 131--139. 1990.

\bibitem[O'Connell(1998)]{oc}
N.~O'Connell.
\newblock Some large deviation results for sparse random graphs.
\newblock \emph{Probab. Th. Rel. Fields}, 110:\penalty0 277--285, 1998.

\bibitem[Palka(1984)]{Palka84}
Z.~Palka.
\newblock On the number of vertices of given degree in a random graph.
\newblock \emph{J. Graph Theory}, 8:\penalty0 167--170, 1984.

\bibitem[Rai\v{c}(2007)]{raic}
M.~Rai\v{c}.
\newblock CLT related large deviation bounds based on Stein's method.
\newblock \emph{Adv. Appl. Prob.}, 39:\penalty0 731--752, 2007.

\end{thebibliography}
\end{document}